\documentstyle[12pt,amstex]{article}
\newcommand{\qed}{\hfill \mbox{$\Box$}}
 \newcommand{\jwsa}{$JW^*$-algebra}
  \newcommand{\pf}{{\it Proof}.}
\newcommand{\tp}[3]{\{#1#2#3\}}
   \newcommand{\tpc}[3]{\{#1,#2,#3\}}

\newtheorem{theorem}{Theorem}[section]

\newtheorem{lemma}[theorem]{Lemma}

\newtheorem{proposition}[theorem]{Proposition}

\date{}

\title{Normal contractive projections preserve type}
\author{Cho-Ho Chu\footnote{Partially supported by Ministerio de Education y Cultura, Spain,
grant PB 98-1371}\,,
 Matthew Neal\footnote{Supported in part by NSF
 grant
 DMS-0101153 \newline \indent
}
\,  and Bernard Russo$^\dag$
}

\begin{document}

\maketitle
\begin{abstract}

Given a JBW*-triple $Z$ and a normal contractive projection $P:Z
\longrightarrow Z$, we show that the (Murray-von Neumann) type of
each summand of $P(Z)$ is dominated by the type of $Z$.
\end{abstract}

\begin{center}
{\bf Introduction}
\end{center}

Contractive projections play a useful role in the theory of
operator algebras and Banach spaces. The ranges of contractive
projections on C*-algebras form an important subclass of those complex
Banach
spaces whose open unit balls are bounded symmetric domains.
 An important feature of these spaces is
that they are equipped with a Jordan {\it triple} product, induced
by the Lie algebra of the automorphism group of the open unit
ball. Known as $JB^*$-triples, they have been shown to be the
appropriate category in which to study contractive projections;
indeed the fact that the category of $JB^*$-triples is stable
under contractive projections played a key role in their structure
theory.

 Recently, contractive projections on von Neumann
algebras have arisen in the study of operator spaces as well as
the theory of harmonic functions on locally compact groups.  In
\cite{NR}, a family of Hilbertian operator spaces were studied and
used to classify, in an appropriate sense, the ranges of
contractive projections on $B(H)$ which are atomic as Banach
spaces. In \cite{C}, it was shown that the Banach space of bounded
matrix-valued harmonic functions on a locally compact group is the
range of a contractive projection on a type I finite von Neumann
algebra. It has also been shown in \cite {CL} that the Banach
space of harmonic functionals on the Fourier algebra of a locally
compact group $G$ is the range of a contractive projection on the
group von Neumann algebra $VN(G)$.

There is a Murray-von Neumann type classification for
$JBW^*$-triples, that is,  $JB^*$-triples which are the dual of a
Banach space. In view of the fact, noted above, that the range of
a contractive projection on a $JB^*$-triple is again a
$JB^*$-triple \cite{FriRus85,K,S}, the above investigations point
to a natural and important question, namely, how is the Murray-von
Neumann classification of the domain affected by a contractive
projection? More precisely, given a $JBW^*$-triple $Z$ of type
$X$, where $X=I,II,\mbox{ or } III$, is the range of a normal
contractive projection on $Z$ of type $Y$ with $Y\le X$, meaning
each summand of the range is of type $\le X$ ? In this paper, we
answer this question affirmatively. We shall see that it suffices
to prove this for $JW^*$-triples, that is, for $JBW^*$-triples
which are linearly isometric to a weak operator closed subspace of
$B(H)$, stable for the triple product $xy^*z+zy^*x$, where $B(H)$
is the von Neumann algebra of bounded operators on a Hilbert space
$H$.

Tomiyama \cite{TO} has analysed the type structure of the range
of a contractive projection which is a von Neumann subalgebra of
the domain. His arguments depend on the crucial fact that the
range is a {\it subalgebra}. In our investigation, the range, which
automatically
has an algebraic structure, need
not be a subalgebra nor even a {\it subtriple}. This adds both generality
and
 complexity to our question.

This paper is organized as follows. Section~\ref{sect:1} is
devoted to background and motivation for the problem. In
section~\ref{sect:2} we consider, as a preliminary tool,
contractive projections on $JW^*$-algebras.
Propositions~\ref{prop:3} and ~\ref{prop:3.1} show that if the
image of a normal contractive projection on a $JW^*$-algebra is a
$JW^*$-subalgebra (not necessarily with the same identity), then
the properties of being semifinite or of type I are passed on from
the domain to the image. In section~\ref{sect:3} we study normal
contractive projections on a von Neumann algebra of type I and
show in Proposition~\ref{prop:1prime}
 that the
image is isometric to a $JW^*$-triple of type I. It is necessary first to
prove this
(in Proposition~\ref{prop:1})
in the special case when the projection is the
Peirce 2-projection with respect to a partial isometry. Our main results,
that
 normal contractive projections on $JW^*$-triples preserve both type I and
semifiniteness,
 appear
in section~\ref{sect:4} as Theorems~\ref{thm:2} and ~\ref{thm:3}. Again,
Propositions~\ref{prop:3.3} and ~\ref{prop:5.2} deal with the special case
of
a Peirce 2-projection. Although Propositions~\ref{prop:3} through
~\ref{prop:5.2} are
each a
special case of Theorem~\ref{thm:2} or ~\ref{thm:3}, they are essential
steps
in the proofs of these theorems and they
 are  new and of interest.
In section~\ref{sect:5}, we extend Theorems~\ref{thm:2} and~\ref{thm:3} to
arbitrary $JBW^*$-triples, and consider the case of atomic
$JBW^*$-triples.

\section{Motivation and Background}\label{sect:1}

Let $M$ be a von Neumann algebra and let $N$ be a von Neumann subalgebra
of $M$ containing the identity element of $M$.  A positive linear map $E:
M\rightarrow N$ satisfying $Ex=x$ for $x\in N$ and $E(axb)=aE(x)b$
for $x\in M$ and $a,b\in N$ is called a {\it conditional expectation}.
Conditional expectations have played some fundamental roles
in the theory of von Neumann algebras, for instance in V.\ Jones'
theory of subfactors. Work in the 1950s of Tomiyama and
Nakamura-Takesaki-Umegaki established that conditional expectations
are idempotent, contractive, and completely positive mappings, and they
preserve type when normal; see the survey
paper of Stormer \cite{Stormer97}.
Conversely (\cite[10.5.85]{KadRin86}), a unital contractive projection
from one
$C^*$-algebra onto a unital $C^*$-subalgebra extends to a normal
conditional expectation
on the universal enveloping von Neumann algebra, and is in particular a
conditional expectation on the $C^*$-algebra.

A type theory for weakly closed
Jordan operator algebras, based on modularity
of the lattice of projections, and
parallel to the type classification theory for von Neumann algebras,
 was introduced and developed in the 1960s
by Topping \cite{T} and Stormer \cite{Stormer66}. In particular, Stormer
showed that a
$JW$-algebra is of type I if and only if its enveloping von Neumann
algebra is of type I.
This was extended to types II and III by Ayupov in 1982 \cite{A}.
In some cases
 the $JW$-algebra in these results is required to be reversible.

A special case of a result of Choi-Effros in 1977 \cite{ChoEff77}, of
fundamental
importance in the rapidly advancing
theory of operator spaces, states that the range of a
unital completely positive projection on a $C^*$-algebra,
while not in general a subalgebra, nevertheless carries the structure of a
$C^*$-algebra. The proof hinges on a conditional expectation formula
(needed to prove
that the abstract product is associative) which is established using the
Kadison-Schwarz inequality for positive linear maps. We note such a
projection  is completely contractive.

A special case of a result of Effros-Stormer in 1979 \cite{EffSto79}
states that the
range of a unital positive projection on a $C^*$-algebra, while not in
general a Jordan
subalgebra, carries a natural Jordan algebra structure. As before, the
proof depends
on a conditional expectation formula (needed to prove that the abstract
product satisfies the Jordan identity), and such a projection is
contractive.

The above results raised the question of what algebraic structure existed
in the range of an arbitrary contractive projection on a $C^*$-algebra.
A special case of a result of Friedman and Russo in 1983 
states
that the range of such a projection is linearly isometric to a subspace,
closed under the triple product $xy^*z+zy^*x$, of the second
dual of the
$C^*$-algebra. Because of the
lack of an order structure and hence the unavailability of the
Kadison-Schwarz
inequality, new techniques were needed and developed by Friedman-Russo
 in their theory of ``operator algebras without order''
(\cite{FriRus83a}),
 including some conditional expectation formulas for the triple product
(\cite[Corollary 1]{FriRus84}).

During the 1980s,  the theory of $JB^*$-triples was developed extensively;
for a summary, see the survey \cite{Russo94}. In
particular, a
 type I theory was developed for $JBW^*$-triples by Horn in his thesis in
1984. In this theory, idempotents (projections) were replaced by
tripotents (which are abstraction of partial isometries), and the
reduced algebra $pAp$ was replaced by the Peirce 2-space of a
tripotent. Of special importance here is the algebraic fact that
such a Peirce 2-space has an abstract structure of a Jordan
algebra, and moreover Horn has proved that a $JBW^*$-triple is of
type I if, and only if, it contains a complete tripotent whose
Peirce 2-space is a Jordan algebra of type I. The remarkable
structure theorem of Horn states that type I $JBW^*$-triples are
isometric to direct sums of tensor products of a commutative von
Neumann algebra by a Cartan factor.

We now recall some definitions. A {\it Jordan triple system} is a
complex vector space $V$ with a {\em Jordan triple product}
$\{\cdot,\cdot,\cdot\} : V \times V \times V \longrightarrow V$
which is symmetric and linear in the outer variables, conjugate
linear in the middle variable and satisfies the Jordan triple
identity
\[\{a,b,\{x,y,z\}\} = \{\{a,b,x\},y,z\} - \{x,\{b,a,y\},z\} +
\{x,y,\{a,b,z\}\}. \]
A complex Banach space $Z$ is called a $JB^*\text{\it -triple}$ if
it is a Jordan triple system such that for each $z\in Z,$ the
linear map $$ z \square z: v\in Z\mapsto \{z,z,v\}\in Z $$ is
Hermitian, that is, $\|e^{it(z \square z)}\| = 1$ for all $t \in
\Bbb R$, with non-negative spectrum and $\Vert z\square z\Vert
=\Vert z\Vert ^2.$ A $JB^*\text{-triple}\; Z$ is called a
$JBW^*\text{\it -triple}$ if it is a dual Banach space, in which
case its predual is unique, denoted by $Z_*,$ and the triple
product is separately weak* continuous.  The second dual $Z^{**}$
of a $JB^*\text{-triple}$ is a $JBW^*\text{-triple.}$ A
norm-closed subspace of a JB*-triple is called a {\it subtriple}
if it is closed with respect to the triple product.   A JBW*-triple is
called
a {\it JW*-triple} if it can be embedded as a subtriple of some $B(H)$.

 The $JB^*\text{-triples}$ form a large class of Banach spaces which
  include $C^*\text{-algebras,}$ Hilbert spaces and spaces of rectangular
matrices. The triple product in a C*-algebra $\cal A$ is given by
$$ \{x,y,z\} = \,\frac 12\; (xy^*z+ zy^*x). $$ In fact, $\cal A$
is a Jordan algebra in the product $$ x \circ y = \frac{1}{2}(xy +
yx)$$ and we have $\{x,y,z\} = (x \circ y^*)\circ z + (y^* \circ
z)\circ x - (z \circ x)\circ y^*.$ A norm-closed subspace of a
$C^*$-algebra is called a JC*-{\it algebra} if it is also closed
with respect to the involution $*$ and the Jordan product $\circ$
given above. A $JC^*\text{-algebra}$ is called a $JW^*\text{-{\it
algebra}}$ if it is a dual Banach space.

An element $e$ in a JB*-triple $Z$ is called a {\it tripotent} if
$\{e, e, e \}=e$ in which case the map $e\square e :  Z
\longrightarrow Z$ has eigenvalues $0,\, {1\over 2}$ and $ 1$, and
we have the following decomposition in terms of eigenspaces $$
Z=Z_2(e)\oplus Z_1(e)\oplus Z_0(e) $$ which is called the {\it
Peirce decomposition} of $Z$. The ${k\over 2}$-eigenspace $Z_k(e)$
is called the {\it Peirce k-space}. The {\sl Peirce projections}
from $Z$ onto the Peirce k-spaces are given by $$ P_2(e) = Q^2(e),
\qquad P_1(e)= 2(e\square e- Q^2(e)), \qquad P_0(e)= I -2e\square
e+Q^2(e) $$ where $Q(e)z = \{ e, z, e\}$ for $z\in Z$. The Peirce
projections are contractive.

In later computation, we will use frequently the {\em Peirce
rules}
\[ \{Z_i (e)\, Z_j (e)\, Z_k (e) \} \subset Z_{i-j+k} (e) \]
where $Z_l (e) = \{0\}$ for $l \neq 0,1,2$. We note that the
Peirce 2-space $Z_2(e)= P_2(e)(Z)$ is a Jordan Banach algebra with
identity $e$, the Jordan product $a\circ b = \{ a, e, b\}$ and
involution $a^{\#} =\{e, a, e\}$ which satisfy $$\|a^{\#}\| =
\|a\|; \qquad \|\{a, a^{\#}, a \} \| = \|a\|^3$$ where $\{x,y,z\}
= (x \circ y^*)\circ z + (y^* \circ z)\circ x - (z \circ x)\circ
y^*$, in other words, $Z_2(e)$ is a unital {\it JB*-algebra}. A
JB*-algebra having a predual is called a {\it JBW*-algebra}. As
shown in \cite{W}, the self-adjoint parts of  JB*-algebras (resp
JBW*-algebras) are exactly the {\it JB-algebras} (resp {\it
JBW-algebras}).
 For definitions and basic results
about JB-algebras, we refer the reader to \cite{HS}. If $Z =
Z_2(e)$, then $e$ is called {\it unitary}. If $Z_0(e) = \{0\}$,
then the tripotent $e$ is called {\it complete}. Two tripotents
$u$ and $v$ are said to be {\it orthogonal} if $u \square v = 0$.
The elements of the predual $Z_*$ of a JBW*-triple $Z$ are exactly
the {\it normal} functionals on $Z$, that is, the continuous
linear functionals on $Z$ which are additive on orthogonal
tripotents.

Given an orthogonal family of tripotents $\{e_i\}_{i \in \Lambda}$
in a JB*-triple $Z$, we can form a {\it joint Peirce
decomposition} $$Z = \bigoplus _{i,j \in \Lambda} Z_{ij}$$ where
Peirce spaces $Z_{ij}$ are defined by $$Z_{ii} = Z_2(e_i)\,, \quad
Z_{ij} = Z_1(e_i) \cap Z_1(e_j) \quad (i\neq j)$$ $$Z_{i0} = Z_1
(e_i) \cap \bigcap_{j \neq i} Z_0 (e_j)\,, \quad Z_{00} =
\bigcap_i Z_0(e_i).$$ We have, for $z_{ij} \in Z_{ij}$ and $e =
\sum e_i$, $$(e_k \square e)(z_{ij}) = (e_k \square e_k)(z_{ij}) =
\left\{
\begin{array}{ll}
0 & \mbox{if $k \notin \{i,j\}$}\\
 \frac{1}{2} z_{ij} & \mbox{ if
$k \in \{i,j\}$}.
\end{array} \right.
$$

JBW*-triples have an abundance of tripotents. In fact, given a
JBW*-triple $Z$ and $f$ in the predual $Z_*$, there is a unique
tripotent $v_f \in Z$, called the {\it support tripotent} of $f$,
such that $f \circ P_2(v_f) = f$ and  the restriction $f|_{
Z_2(v_f)}$ is a {\it faithful positive} normal functional.

The Murray-von Neumann classification of the von Neumann algebras
can be extended to that of JBW*-triples  and, a JBW*-triple  can
be decomposed into a direct sum of type $j$ ($j =I, II, III$)
summands (see \cite{H1,H2}). A JBW*-triple is called {\it
continuous} if it does not contain a type I summand in which case,
it is a direct sum of a JW*-algebra $H(A, \alpha)$ and a weak*
closed right ideal of a continuous von Neumann algebra, as shown
in \cite{H2}, where $$H(A, \alpha) = \{ a \in A: \alpha (a) = a
\}$$ is the fixed-point set of a period 2 weak* continuous
antiautomorphism $\alpha$ of a von Neumann algebra $A$. It follows
that continuous JBW*-triples are JW*-triples.

 A JBW*-triple $Z$ is
called {\it type} I if it contains an abelian tripotent $e$ such
that $Z= U(e)$ where $U(e)$ denotes the weak* closed triple ideal
generated by $e$. We recall that a tripotent $e$ is said to be
{\it abelian} if the Peirce 2-space $P_2(e)(Z)$ is an abelian
triple which is equivalent to saying that $P_2(e)(Z)$ is an
associative JBW*-algebra in the usual Jordan product $x\circ y =
\{x,e,y\}$. Horn \cite[4.14]{H3} has shown that a JBW*-triple is
type I if, and only if, every weak*-closed triple ideal
 contains an abelian tripotent.

An important class of type I JBW*-triples are the following six
types of {\it Cartan factors}:

\begin{tabular}{rl}
type 1 & $B(H,K)$~~  {\em with triple product}~ $\{x,y,z\}
=\frac{1}{2}(xy^*z +zy^*x),$\\ type 2 & $\{z\in B(H,H): z^t =
-z\},$
\\ type 3 &$\{z\in B(H,H): z^t = z\},$ \\ type 4 &spin factor, \\
type 5 &$M_{1,2}({\cal O})$~~  {\em with triple product}~
$\{x,y,z\} = \frac{1}{2}(x(y^*z) + z(y^*x)),$\\ type 6 &$M_3({\cal
O})$
\end{tabular}\\

\noindent where $B(H,K)$ is the Banach space of bounded linear
operators between complex Hilbert spaces $H$ and $K$, and $z^t$ is
the transpose of $z$ induced by a conjugation on $H$. Cartan
factors of type 2 and 3 are subtriples of $B(H,H)$, the latter
notation is shortened to $B(H)$. The type 3 and 4 are Jordan
algebras with the usual Jordan product $x \circ y = \frac{1}{2}(xy
+ yx)$. A {\em spin factor\/} is a Banach space that is equipped
with a complete inner product $\langle \cdot, \cdot \rangle$ and a
conjugation $j$ on the resulting Hilbert space, with triple
product $$ \{ x, y, z\} = \frac{1}{2} (\langle x,y\rangle z +
\langle z,y \rangle x - \langle  x, jz \rangle jy)$$ such that the
given norm and the Hilbert space norm are equivalent.

 By Horn's
result in
\cite{H1}, a JBW*-triple
$Z$ is of type I if, and only if, it is linearly isometric to an
$\ell^{\infty}$-sum
$\bigoplus_{\alpha} L^{\infty}(\Omega_{\alpha}) \otimes
C_{\alpha}$ where $C_{\alpha}$ is a Cartan factor. Such a type I
JBW*-triple will be called {\it type} I$_{fin}$ if each Cartan
factor
$C_{\alpha}$ is finite-dimensional.
 It has been
shown in \cite{CM} that a JBW*-triple  $Z$ is type I$_{fin}$ if,
and only if, its predual $Z_*$ has the Dunford-Pettis property. We
recall that a Banach space $W$ has the {\it Dunford-Pettis
property} if every weakly compact operator on $W$ is completely
continuous. Such property is inherited by complemented subspaces.

Horn's type I structural result above also shows that a
JBW*-algebra is type I as a JBW*-triple if and only if its
self-adjoint part is a type I JBW-algebra in the sense of
\cite{HS}.

\begin{lemma}\label{sub}
Let $Z$ be a JBW*-subtriple of a type I$_{fin}$ JBW*-triple.
Then $Z$ is type I$_{fin}$.
\end{lemma}
\pf\ By \cite[Corollary 6]{CM}, $Z_*$ has the Dunford-Pettis
property.\qed

\medskip

We will begin our investigation of contractive projections in the
next section.  A contractive projection $P: Z \longrightarrow Z$
on a JB*-triple $Z$ is a bounded linear map such that $P^2 = P$
and $\|P\| \leq 1$. We will exclude the trivial case of $P=0$
which then implies $\|P\| =1$. Given such a contractive projection
$P$ on $Z$ with triple product $\{\cdot, \cdot, \cdot,\}$, one can
show, using the holomorphic characterization of JB*-triples
\cite{K,S}, that the range $P(Z)$ is also a JB*-triple in the
triple product $$ [x,y,z] = P\{x,y,z\} \qquad (x,y,z \in P(Z)).$$
Moreover, one has the following conditional expectation formula:
$$ P\{Px, Py, Pz\} = P\{Px, y, Pz\} \qquad (x,y,z \in Z).$$ The
above result has also been proved in \cite{FriRus85} for
subtriples of C*-algebras, via an operator algebra approach which
also yields the formula $$ P\{Px, Py, Pz\} = P\{x, Py, Pz\}.$$ A
weak* continuous projection on a JBW*-triple is called {\it
normal}.

\section{Contractive projections on $JW^*$-algebras}\label{sect:2}

In this section, we consider a JW*-algebra $A \subset B(H)$ with
positive part $A^+$, inheriting various topologies of $B(H)$. A
positive linear functional $\varphi$ of $A$ is called a {\it trace }
if $\varphi
(sxs) = \varphi (x)$ for all symmetries $s \in A$ and all $x \in
A^+$, where a {\it symmetry} in $A$ is a self-adjoint element $s$
such that $s^2$ is the identity in $A$. By \cite{A}, every normal
trace on $A$ can be extended to a normal trace on its enveloping
von Neumann algebra. Further, if $\varphi$ is faithful, so is its
extension. In the sequel, our $JW^*$-subalgebras need not
 have the same identity element as the $JW^*$-algebras which
contain them.

The following lemma is a special case of Lemma \ref{sub}, but the
proof below is {\it intrinsic} without using the Dunford-Pettis
property.

\begin{lemma}\label{lem:2.3}
Every JW*-subalgebra of a  type I$_{fin}$
JW*-algebra is
of type I$_{fin}$.
\end{lemma}
\pf\ Let $A$ be a JW*-subalgebra of a type I$_{fin}$ JW*-algebra
$B$. Then $A$ is finite since it is a subalgebra of a finite
algebra. Let $p \in A$ be a projection. Then $pAp$ is a subalgebra
of the type I$_{fin}$ algebra $pBp$. Suppose, for contradiction,
that $pAp$ contains no abelian projection. By cutting down to a
homogeneous summand, we may assume that $pBp$ is homogeneous. Then
by \cite[Theorem 17]{T}, $p$ can be decomposed into any number of
mutually orthogonal and strongly equivalent projections in $pAp$.
Since equivalent projections in $pAp$ are also equivalent in
$pBp$, and since in a homogeneous type I$_{fin}$ algebra, there
are at most a fixed number of mutually orthogonal and  strongly
equivalent projcetions, we have a contradiction. So $pAp$ contains
an abelian projection and $A$ is type I$_{fin}$.\qed

\begin{lemma}\label{id}
Let $(A,\circ)$ be a JW*-algebra with identity $1$ and let $P: A
\longrightarrow A$ be a contractive projection. If $P(A)$ contains
a unitary tripotent $u$ of $P(A)$, then $P1 = P(uu^*u^*u) =
P(u\circ u^*)$. In addition, if $u$ is a projection in $A$, then
$P1 = u$.
\end{lemma}
\pf\ Recall that the triple product in $P(A)$ is given by
$$[x,y,z] = P\{x,y,z\}.$$ Since $u$ is a unitary tripotent in
$P(A)$, we have by the main identity
\begin{eqnarray*}
P1 & = & [P1,u,[u,u,u]]\\
&=& [[P1,u,u],u,u]-[u,[u,P1,u],u]+[u,u,[P1,u,u]]\\
&=&P1-[u,[u,P1,u],u]+P1
\end{eqnarray*}
and by the conditional expectation formula,
\begin{eqnarray*}
P1 & = & [u,[u,P1,u],u]\\
&=& P\{u,P\{u,P1,u\},u\}\\
&=& P\{u,P(u^2),u\} = P\{u,u^2,u\} = P(uu^*u^*u).
\end{eqnarray*}
Also, $P1 =[u,u,P1] = P\{u,u,P1\}= P\{u,u,1\} = P(u\circ
u^*)$.\qed

\medskip

\noindent {\bf Remark.}~ The above result shows that there is
at most one unitary tripotent in $P(A)$ which is a projection
in $A$. If $P(A)$ is a JW*-subalgebra of $A$, then the identity
in $P(A)$ is a projection in $A$ and $P(1_A)=1_{P(A)}$.\\

As in \cite{T}, $A$ is said to be {\it modular} if its projections
form a modular lattice in which case $A$ admits a centre-valued
trace, and therefore a separating family of normal traces. It has
been shown in \cite{A} that a JW*-algebra is modular if, and only
if, its enveloping von Neumann algebra is finite. For this reason,
we propose from now on to replace the term ``modular'' by the more
common term ``finite'' throughout. A projection $p$ in a
JW*-algebra $A$ is called {\it finite} if the JW*-algebra $pAp$ is
finite.\\

We recall that, for a net $(x_{\alpha})$ in a von Neumann
algebra, we have
$$x_{\alpha} \longrightarrow 0~ strongly \Leftrightarrow
x_{\alpha}^*x_{\alpha} \longrightarrow 0~ weakly$$
$$x_{\alpha} \longrightarrow 0~ strongly^* \Leftrightarrow
x_{\alpha}^*x_{\alpha} +x_{\alpha}x_{\alpha}^* \longrightarrow
0~ weakly.$$
Plainly, strong* convergence implies strong convergence.

\begin{lemma}
Let $A$ be a JW*-algebra. The following conditions are equivalent:
\begin{description}
\item (i)  $A$ is finite.
\item (ii) The map $x \in A \mapsto x^* \in A$ is strongly
continuous on bounded spheres in the  enveloping von Neumann
algebra of $A$.
\end{description}
\end{lemma}
\pf\ $(i)\Rightarrow (ii)$. Let $\cal F$ be a separating family of
normal traces of $A$. Then the family $\tilde {\cal F} =
\{\tilde{\varphi}: \varphi \in {\cal F}\}$ of normal tracial
extensions of traces in $\cal F$ is separating on the enveloping
von Neumann algebra $\cal A$ of $A$ (cf. \cite[proof of Theorem
3]{A}). We have $x_{\alpha} \rightarrow 0$ strongly $\Rightarrow$
$x_{\alpha}^*x_{\alpha} \rightarrow 0$ weakly $\Rightarrow$
$\tilde{\varphi}(x_{\alpha}x_{\alpha}^*) =
\tilde{\varphi}(x_{\alpha}^*x_{\alpha}) \rightarrow 0$ for all
$\tilde{\varphi} \in \tilde{\cal F}$. Hence
$x_{\alpha}x_{\alpha}^* \rightarrow 0$ weakly, that is,
$x_{\alpha}^* \rightarrow 0$ strongly.\\

$(ii)\Rightarrow (i)$. If $A$ is not finite, then by \cite[Lemma
23]{T}, there is an infinite orthogonal sequence $\{p_n\}$ of
projections in $A$ such that, for every $n$, $$p_1 = s_n p_n s_n$$
where $s_n$ is a symmetry. Given any normal state $\psi$ of the
enveloping von Neumann algebra of $A$, we have $$\sum \psi (p_n) =
\psi (\sum p_n) \leq \psi (1) < \infty.$$ So $\psi
((s_np_n)^*(s_np_n)) = \psi (p_n) \rightarrow 0.$ But $\psi((s_n
p_n)(s_np_n)^*) = \psi (p_1) \not \rightarrow 0$. So the map $x
\mapsto x^*$ is not strongly continuous on the unit ball.\qed

\begin{lemma} \label{lem:3.4}Let $P: A \longrightarrow A$ be a
 contractive projection on a JW*-algebra $A$
such that
$P(A)$ is a JW*-subalgebra of $A$. Then
\begin{description}
\item (i) $P(x\circ x^*) \geq 0$ for all $x \in A$;
\item (ii) $P(a\circ x) = P(a)\circ x$ for $a \in A, x \in
P(A)$;
\item (iii) If $P$ is normal, then $P$ is strongly* continuous
on bounded spheres.
\end{description}
\end{lemma}
\pf\ (i) By Lemma \ref{id}, $P1$ is the identity in $P(A)$. Let
$\varphi$ be a state of $P(A)$. Then $\varphi\circ P(1) = \varphi
(P1) = 1$ implies that $\varphi \circ P$ is a state of $A$. Hence
$\varphi (P(x \circ x^*)) \geq 0.$ As $\varphi$ was arbitrary, we
have $P(x\circ x^*) \geq 0$. This implies that $P$ is
self-adjoint.\\

(ii) This is proved in \cite{EffSto79}. We give a short
alternative proof here. We have $P(a\circ x) = P(a\circ(x\circ
P1)) = \frac{1}{2}P\{a,x,P1\} + \frac{1}{2}P\{a,P1,x\} =
\frac{1}{2}(P\{Pa,x,P1\} + P\{Pa, P1, x\})
=\frac{1}{2}(\{Pa,x,P1\} + \{Pa, P1, x\}) = Pa\circ x$.\\

(iii) Let $x_{\alpha} \rightarrow 0$ strongly*. Then
$x_{\alpha}\circ x_{\alpha}^* \rightarrow 0$ weakly and hence
$P(x_{\alpha}\circ x_{\alpha}^*) \rightarrow 0$ weakly. Using
(i), (ii), and the self-adjointness of $P$, we have
$$0 \leq P((Px_{\alpha} - x_{\alpha})\circ (Px_{\alpha} -
x_{\alpha})^*) = P(x_{\alpha}\circ x_{\alpha}^*) -
P(x_{\alpha})\circ P(x_{\alpha})^*$$
which implies that $P(x_{\alpha})\circ P(x_{\alpha})^*
\rightarrow 0$ weakly. \qed

\medskip

A JW*-algbra $A$ is {\it semifinite} if every nonzero projection in
$A$ contains a nonzero finite projection. This is equivalent to
saying that $A$ does not contain any type III summand.

\begin{proposition}\label{prop:3}
 Let $A$ be a semifinite (resp. finite) JW*-algebra and $P$
a normal contractive projection on $A$ such that
$P(A)$ is a JW*-subalgebra of $A$. Then
$P(A)$ is semifinite (resp. finite).
\end{proposition}
\pf\ Suppose $P(A)$ is of type III. We show that $P(A) = 0$. Let
$e\in A$ be a finite projection. Suppose $P(e) \neq 0$. We have
$P(e)\ge 0$ and by spectral theory there
 exists a nonzero projection $p \in P(A)$ such that
$\lambda p \leq P(e)$ for some $\lambda > 0$. Let
$(x_{\alpha})$ be a bounded net in $pP(A)p$ converging to $0$
strongly, in the enveloping von Neumann algebra $\cal A$ of
$A$. Since $e$ is also a finite projection in $\cal A$
(cf.\cite[Corollary 3.2]{A}), by
\cite[p.97-98]{SA}, the nets $(x_{\alpha}^*e)$ and
$(ex_{\alpha}^*)$ converge to $0$ strongly in $\cal A$.
Since the nets $(ex_{\alpha})$ and $(x_{\alpha}e)$ both
converge to 0 strongly, we have $(x^*_{\alpha}e)$ and
$(ex^*_{\alpha})$ both converging to 0 strongly*. Therefore
$\{e,x_{\alpha},P(e)\}
\longrightarrow 0$ strongly* in $A$.
 Since $P$ is strongly* continuous on bounded spheres
and $P(A)$ is in particular a subtriple of $A$, we have
$$P(e)x^*_{\alpha}P(e) =P\tpc{Pe}{x_\alpha}{P(e)}=
 P\{e,x_{\alpha},P(e)\}
\longrightarrow  0$$ strongly*, and hence strongly. It follows
that
$$ [pP(e)p +(1-p)]x^*_{\alpha}[pP(e)p +(1-p)] =
pP(e)x^*_{\alpha}P(e)p \longrightarrow 0$$
strongly which gives
$$x^*_{\alpha} =
[pP(e)p+(1-p)]^{-1}pP(e)x^*_{\alpha}P(e)p[pP(e)p+(1-p)]^{-1}
\longrightarrow 0$$
strongly, implying that $pP(A)p$ is finite and contradicting
that  $P(A)$ is type III.  Hence $P$ vanishes on every finite
projection in $A$ and $P(A) = 0$.

If we apply the above argument to the identity element of a finite
$A$, we obtain that $P(A)$ is finite. \qed

\medskip

It will follow from Theorems~\ref{thm:2} and ~\ref{thm:3} in
section~\ref{sect:4}
that Proposition~\ref{prop:3}, and
Proposition~\ref{prop:3.1} which follows, remain
 true without
 the assumption that $P(A)$ is a subalgebra.
The proof of Proposition~\ref{prop:3.1} is an adaptation to the
Jordan algebra setting of the proof for von Neumann algebras in
\cite{TO}.

\begin{proposition}\label{prop:3.1}
Let $P$ be a normal contractive projection on a type I \jwsa\ $A$
and suppose $P(A)$ is a $JW^*$-subalgebra of $A$.  Then $P(A)$ is
of type I.
\end{proposition}
\pf\
By Proposition~\ref{prop:3}, $P(A)$ is a semifinite
$JW^*$-algebra. Suppose that $P(A)$ contains
a type II summand.
By following $P$ by the projection onto the type II
part, we can assume that $P(A)$ is of type II. We show $P(A) = 0$.
It suffices to show that for any finite projection $q$ in $P(A)$,
we have $qP(A)q=0$. By following $P$ with the projection $q\cdot q$, we
may further assume that $P(A)$ is of type $II_1$. Suppose $B = P(A)
\neq 0$, we deduce a contradiction.

By \cite[Theorems 2 and 5]{A} there are faithful normal semifinite
traces $\tau,\tau_0,\tilde{\tau}_0$
on $B$, $A$, $\tilde{A}$ respectively, where
$\tilde{A}$ is the von Neumann algebra
generated by $A$, such that $\tau$ is finite and $\tilde{\tau}_0$ is
an extension of $\tau_0$.  Since $\tau\circ P$ is a normal positive
functional on $A$, by the Radon-Nikodym theorem \cite[Theorem 2.4]{A1},
there is an
operator
$h\in L^1(A,\tau_0)^+$ such that
\begin{equation}\label{eq:RN}
\tau\circ P(x)=
\tau_0(\tp{h^{1/2}}{x}{h^{1/2}})\mbox{ for }x\in A.
\end{equation}

Note that for self-adjoint $x\in A$ and $y\in B$, $\tau\circ
P(y^2\circ x)= \tau(P(y^2\circ x))=\tau(y^2\circ Px)$ by
Lemma~\ref{lem:3.4}(ii). On the other hand, $\tau\circ
P(\tp{y}{x}{y})=\tau(P\tp{y}{x}{y})=\tau(\tpc{y}{Px}{y})
=\tau(y^2\circ Px)$, the latter by \cite{PS}. Hence $\tau\circ
P(\tp{y}{x}{y})=\tau\circ P(y^2\circ Px)$, for self-adjoint $x\in
A,y\in B$. Applying this to a projection $p\in B$ and using the
extension property and (\ref{eq:RN}), we have
\[
\tilde{\tau}_0\left(h\left[\frac{xp+px}{2}-pxp\right]\right)=0\mbox{
for every }x\in A.
\]

Expanding
 $\tilde\tau_0(h\circ(p\circ x))=\tilde\tau_0(h\circ(pxp))$ and using
the associative trace properties of $\tilde\tau_0$ yields
\[
\tilde\tau_0\left(\frac{xhp+xph+xhp+xph}{4}\right)=\tilde\tau_0\left(
\frac{xphp+xphp}{2}\right).
\]
Hence,
\[
\tilde{\tau}_0(x(ph+hp-2php))=0,
\]
which is the same as $\tau_0(x\circ (ph+hp-2php))=0$. Since this is true
for
all $x\in A$, we have $ph=php=hp$, so that $h$ is affiliated with $B'$,
the commutant of $B$ (see \cite{A1}). Note that, since $p$ is a finite
projection,
all of the strong products above are in $L^1(\tilde{A},\tilde\tau_0)$.

Since $h\in L^1(A,\tau_0)$, we may pick a nonzero finite
projection $e\in A\cap B'$
(a spectral projection of $h$).
It is easy to see that $eB$ is a $JW$-subalgebra of $A$ and that $eB''
=(eB)''$.  By \cite[Theorem 8]{A}, $B$ of type $II_1$ $\Rightarrow B''
\mbox{ of type }II_1\Rightarrow eB''\mbox{ of type }II_1\Rightarrow
(eB)''\mbox{ of type }II_1\Rightarrow eB\mbox{ of type }II_1$, the latter
since $eB$ is reversible.  But $eB=eBe\subset eAe$ is of type $I_{fin}$ by
Lemma~\ref{lem:2.3}, giving a contradiction. Hence  $B=0$.
\qed

\section{Contractive projections on von Neumann algebras}\label{sect:3}

\begin{proposition} \label{prop:1}
Let $M$ be a von Neumann algebra of type I and let $e$ be a partial
isometry of $M$. Then the Peirce 2-space $P_2(e)M$ is a \jwsa\ of type I.
\end{proposition}
\pf\
We note that $P_2(e)M$ is a  von Neumann algebra with identity $e$
under the product $x\cdot
y=xe^*y$ and involution $x^\sharp=ex^*e$ as well as a
\jwsa\ under $x\circ y=\tp{x}{e}{y}=
(xe^*y+ye^*x)/2$ and $x^\sharp$. Also, $(P_2(e)M,\cdot)$ is a von
Neumann algebra of type I if and only if $(P_2(e)M,\circ)$
is a \jwsa\ of type I.

Now suppose that $v$ is a nonzero central projection
in $(P_2(e)M,\cdot)$. Below we shall verify the following:
\begin{description}
\item (i) $v$ is a tripotent in $M$.
\item (ii) $v\cdot P_2(e)M=P_2(v)M$ as sets.
\item (iii) The identity map $:(v\cdot P_2(e)M,\cdot)\rightarrow
(P_2(v)M,\times)$, where $x\times y=xv^*y$ and $x\mapsto vx^*v$ is the
involution in $(P_2(v)M,\times)$,
is a $^*$-isomorphism of von Neumann algebras.
\item (iv) $(P_2(v)M,\times)$ has a non-zero abelian projection.
\end{description}

Assuming that (i)-(iv) have been proved, if
 there is a nonzero central projection $v$ such that $v\cdot P_2(e)M$ is a
continuous
von Neumann algebra,
we obtain a contradiction that it contains a nonzero abelian projection.
So
$(P_2(e)M,\cdot)$ is  a von Neumann algebra of type I.

It remains to verify (i)-(iv) above.

Since $v=v\cdot v=v^\sharp$, we have $v=v\cdot v\cdot v=ve^*ve^*v
=v(ev^*e)^*v=vv^*v$. This proves (i).

Since $v=ee^*ve^*e$, we have $v=ee^*v=ve^*e$ so that
$vv^*e=v(ev^*e)^*e=ve^*ve^*e=ve^*v=
v$ and similarly
$ev^*v=v$.  Hence, for $y\in M$,   $v\cdot P_2(e)y=v\cdot P_2(e)y\cdot
v=ve^*(P_2(e)y)e^*v=vv^*ee^*(P_2(e)y)e^*ev^*v= P_2(v)y$ (since $v^*v\le
e^*e$ and
$vv^*\le ee^*$). This proves (ii).

For $x,y\in P_2(v)M$, we have, by (ii), $x=xe^*v$ and $y=ve^*y$.
  Therefore
$ x \cdot y = xe^*y=xe^*ve^*ve^*y=xv^*ve^*y=xv^*y = x \times y$.  As for
the
involution, $ex^*e=e(ve^*xe^*v)^*e=ev^*ex^*ev^*e=vx^*v$. This proves
(iii).

Let $p=vv^*$. We shall show that
\begin{description}
\item (a) There is a non-zero abelian projection $h\in M$ with $h\le p$
and
$c(h)=c(p)$, where $c(\cdot)$ denotes central support.
\item (b) With $z=hv$, $z$ is a non-zero abelian projection in the von
Neumann algebra $(P_2(v)M,\times)$. This will prove (iv).
\end{description}

Since $pMp$ is of type I, there is a non-zero abelian projection
$h\in pMp$ such that $c_{pMp}(h)=p$. Since $h\le p$, we have
$c(h)\le c(p)$. To show equality here, take a central projection
$r\in M$ with $ h \leq r$. Then $pr$ is a central projection in
$pMp$, so that $pr=prp\ge c_{pMp}(h)=p$ and thus $p\le r$ gives
$c(p)\le r$. Taking $r=c(h)$, we get $c(p)\le c(h)$. This proves
(a).

We next show that $z$ is a non-zero projection in $(P_2(v)M,\times)$.
We have $z=hv=phvv^*v=vv^*hvv^*v\in P_2(v)M$,
$vz^*v=v(vv^*hvv^*v)^*v=vv^*hv=hv=z$, $z\times z=
zv^*z=hvv^*hv=hphv=hv=z$, and $zz^*=hvv^*h=h\ne 0$.

It remains to show that for $x,y\in M$, we have
\begin{eqnarray}\label{eq:abelian}
\lefteqn{[z\times (vv^*xv^*v)\times z]v^*[z\times (vv^*yv^*v)\times
z]=}\\\nonumber
&=&[z\times (vv^*yv^*v)\times z]v^*[z\times (vv^*xv^*v)\times
z].
\end{eqnarray}
The left and right sides of (\ref{eq:abelian}) collapse to
$hxv^*hyv^*hv$ and $hyv^*hxv^*hv$ respectively, which are equal
since  $hMh$ is an abelian subalgebra of $M$. For example, the left side
is equal to $$zv^*vv^*xv^*vv^*zv^*zv^*vv^*yv^*vv^*z=
hvv^*xv^*hvv^*hvv^*yv^*hv=hxv^*hyv^*hv. $$
This proves (b), hence (iv) and the Proposition.\qed

\medskip

Let $P$ be a normal contractive projection on a $JBW^*$-triple $Z$
and let $f\in P_*(Z_*)$ have the support tripotent (partial
isometry in this case) $v_f$. Let $P_k=P_k(v_f)$ denote the Peirce
projections induced by $v_f$.
 The following commutativity formulas were proved in \cite{FriRus87}.
These will be used freely in the remainder of the paper.
\begin{itemize}
\item $P_2P=P_2PP_2\,,\quad PP_2=PP_2P\,;$
\item $PP_0=P_0PP_0=PP_0P\,;$
\item $PP_1=PP_1P\,,\quad P_1PP_0=0\,.$
\end{itemize}

In the next lemma, we shall use these formulas to extend the first two
of them to the case where the tripotent is not assumed to be
the support of a
normal functional.
We shall use the fact that,
by Zorn's lemma, every tripotent in a $JBW^*$-triple
$Z$ is the sum of an orthogonal family of tripotents which
are support tripotents of normal functionals on $Z$.

The following lemma is needed in the next section. In this
section, it will be used only in the case that $Z$ is a von
Neumann algebra, considered as a $JW^*$-triple under
$\frac{1}{2}(xy^*z+zy^*x)$.

\begin{lemma}\label{lem:4.1}
Let $P$ be a normal contractive projection on a $JBW^*$-triple $Z$
and suppose that $v$ is a tripotent of the $JBW^*$-triple $P(Z)$.
Choose a set $S=\{f_i:i\in I\}$ of pairwise orthogonal normal
functionals on $P(Z)$ such that $v=\sum_{i\in I}v_{f_i}$, where
$v_{f_i}$ is the support tripotent of $f_i$ in $P(Z)$. Let $w_i$
be the support tripotent of $f_i$ in $Z$, necessarily pairwise
orthogonal, and let $w$ be the partial isometry $\sum_{i\in I}
w_i$. Then
\[
P_2(w)P=P_2(w)PP_2(w)\,,\quad PP_2(w)=PP_2(w)P.
\]
\end{lemma}
\pf\
Since $w=\sum w_i$, we have $P_2(w)=\sum_i P_2(w_i)+\sum_{j\ne
k}P_1(w_j)P_1(w_k)$ and
therefore
\begin{eqnarray*}
\lefteqn{P_2(w)PP_2(w)}\\&=&\sum_{i,i'}P_2(w_i)PP_2(w_{i'})+\sum_{j\ne
k,j'\ne k'}
P_1(w_j)P_1(w_k)PP_1(w_{j'})P_1(w_{k'})\\
&+&\sum_{j'\ne k',\mbox{ all }i}P_2(w_i)PP_1(w_{j'})P_1(w_{k'})+
\sum_{j\ne k,\mbox{ all }i'}P_1(w_j)P_1(w_k)PP_2(w_{i'}).
\end{eqnarray*}

Because $P_2(w_i)P=P_2(w_i)PP_2(w_i)$, by properties of the
joint Peirce decomposition, the first sum reduces to
$\sum_i P_2(w_i)P$ and each term in the third sum is zero.

Each term in the fourth sum is zero as well. Indeed, since in the
following
we may assume $k\ne i'$,
\begin{eqnarray*}
\lefteqn{P_1(w_j)P_1(w_k)PP_2(w_{i'})}\\
&=&P_1(w_j)[I-P_2(w_k)-P_0(w_k)]PP_2(w_{i'})\\
&=&P_1(w_j)[PP_2(w_{i'})-P_2(w_k)PP_2(w_{i'})-P_0(w_k)PP_2(w_{i'})]\\
&=&P_1(w_j)[PP_2(w_{i'})-0-P_0(w_k)PP_0(w_k)P_2(w_{i'})]\\
&=&P_1(w_j)[PP_2(w_{i'})-PP_0(w_k)P_2(w_{i'})]\\
&=&P_1(w_j)[PP_2(w_{i'})-PP_2(w_{i'})]=0.
\end{eqnarray*}

The second sum reduces to $\sum_{j\ne k} P_1(w_j)P_1(w_k)P$.
Indeed, if $k\not\in\{j',k'\}$, then
$[P_1(w_k)PP_1(w_{j'})]P_1(w_{k'}) =0$ since
$P_1(w_{j'})P_1(w_{k'})Z\subset P_0(w_k)Z$ and\\
$P_1(w_k)PP_0(w_k)=0$. Thus the second sum is reduced to\\
$\sum_{j\ne k} P_1(w_j)P_1(w_k)PP_1(w_j)P_1(w_k)$. However,
\begin{eqnarray*}
\lefteqn{P_1(w_j)P_1(w_k)PP_1(w_j)P_1(w_k)}\\
&=&P_1(w_j)P_1(w_k)P[
I-P_2(w_j)-P_0(w_j)]P_1(w_k)\\
&=&P_1(w_j)P_1(w_k)PP_1(w_k)-P_1(w_j)P_1(w_k)PP_2(w_j)P_1(w_k)\\
&&-P_1(w_j)P_1(w_k)PP_0(w_j)P_1(w_k)\\
&=&P_1(w_j)P_1(w_k)P[I-P_2(w_k)-P_0(w_k)]+0+0\\
&=&P_1(w_j)P_1(w_k)P-P_1(w_j)P_1(w_k)PP_2(w_k)-P_1(w_j)P_1(w_k)PP_0(w_k)\\
&=&P_1(w_j)P_1(w_k)P-P_1(w_k)P_1(w_j)PP_2(w_k)\\
&=& P_1(w_j)P_1(w_k)P.
\end{eqnarray*}
This proves the first formula.

For the second formula, we have
\begin{eqnarray*}
PP_2(w)P&=&P\left(\sum_i P_2(w_i)+\sum_{j\ne k}P_1(w_j)P_1(w_k)\right)P\\
&=&\sum_i PP_2(w_i)P+\sum_{j\ne k}PP_1(w_j)P_1(w_k)P\\
&=&\sum_i PP_2(w_i)+\sum_{j\ne k}PP_1(w_j)P_1(w_k)=PP_2(w)
\end{eqnarray*}
since \begin{eqnarray*}
PP_1(w_j)P_1(w_k)P&=&(PP_1(w_j)P)P_1(w_k)P\\
&=&PP_1(w_j)(PP_1(w_k)\\
&=&PP_1(w_j)P_1(w_k).\qed
\end{eqnarray*}

The following lemma is probably known. We include a proof for
completeness.

\begin{lemma}\label{lem:4.2}
Let $p$ be a projection in a $JB^*$-algebra $A$ and let $A_1(p)$
be the Peirce 1-space. Then $A_1(p)\cap A^+=0$
\end{lemma}
\pf\ If $x\in A_1(p)\cap A^+$, then let $y = x^{\frac{1}{2}} \in
A_{sa}$ and let $y=y_2+y_1+y_0$ be its Peirce decomposition with
respect to $p$. Then $x=y_2^2+y_1^2+y_0^2+2(y_2+y_0)y_1$. Since
$x\in A_1(p)$, we have $y_2^2+y_1^2+y_0^2=0$ and because the
JB-algebra $A_{sa}$ is formally real, $y=0$.\qed

\begin{lemma}\label{lem:4.3}
Let $P,Z,v,S,w$ be as in Lemma~{\rm \ref{lem:4.1}}. Then
\begin{description}
\item (a) The map
$Q=P_2(w)P:Z_2(w)\rightarrow Z_2(w)$ is a normal faithful
unital contractive projection with range $P_2(w)P(Z)$.
\item (b) The map $P_2(w)$
is a linear surjective isometry of $P(Z)_2(v)$ onto\\ $P_2(w)P(Z)$
\end{description}
\end{lemma}
\pf\
(a)~ By Lemma~\ref{lem:4.1}, $Q^2=P_2(w)PP_2(w)P=P_2P=Q$ and
$Q(Z_2(w))=P_2(w)PP_2(w)(Z)=P_2(w)P(Z)$. To show that $Q$ is unital, note
first that by
\cite[Lemma 2.7]{FriRus83}, $v_i=w_i+P_0(w_i)v_i$ so that
$w_i \perp (v_i-w_i)$. By taking sums and limits,
 one obtains $(v-w)\perp w$ and $\|v-w\|\le 1$.
Indeed, it is easy to see that for any finite set $F$ of indices,
$\sum_F w_i$ is the support tripotent of the normal functional
$\sum_Ff_i$. Hence,  $\sum_Fw_i \perp \sum_F(v_i-w_i)$ so that
$\sum_F(v_i-w_i)\in Z_0(\sum_Fw_i)$ and
$\|\sum_Fw_i\pm\sum_F(v_i-w_i)\|=1$. By passing to the limit and
noting that each $f_i$ has the value 1 on $w\pm(v-w)$, we have
$\|w\pm(v-w)\|=1$, and since $P_2(w)$ is contractive, $\|w\pm
P_2(w)(v-w)\|\le 1$, and since $w$ is an extreme point,
$P_2(w)(v-w)=0$, that is, $P_2(w)v=w$. Now $v=w+P_1(w)v+P_0(w)v$,
so by \cite[Lemma 1.6]{FriRus85a}, $P_1(w)v=0$ and thus
$v=w+P_0(w)v$ and $v=Pv=Pw+PP_0(w)x$ so that
$Qw=P_2(w)Pw=P_2(w)(v-PP_0(w)v)=P_2(w)v=w$ and $Q$ is unital.

Finally we show that $Q$ is faithful. Suppose
that $b\in Z$, $P_2(w)b\ge 0$, and $P_2(w)Pb=0$  We shall show that
$P_2(w)b=0$. In the
first place, since $P_2(w_i)$ is a positive
 operator on the $JB^*$-algebra $P_2(w)Z$ (\cite[3.3.6]{HS}),
$P_2(w_i)b=P_2(w_i)P_2(w)b\ge 0$ for every $i\in I$.  Since
$P_1(w_k)P_1(w_l)b\perp w_i$, we have
\begin{eqnarray*}
0&=&\langle P_2(w)Pb,f_i\rangle=\langle PP_2(w)Pb,f_i\rangle
=\langle PP_2(w)b,f_i\rangle=\langle P_2(w)b,f_i\rangle\\
&=&\sum_j\langle P_2(w_j)b,f_i\rangle+ \sum_{k\ne l}\langle
P_1(w_k)P_1(w_l)b,f_i\rangle=\langle P_2(w_i)b, f_i\rangle.
\end{eqnarray*}
 Hence $P_2(w_i)b=0$ for all $i$.   Therefore
 $P_2(w)b=\sum_i P_2(w_i)b+\sum_{j\ne k}P_1(w_j)P_1(w_k)b=
\sum_{j\ne k}P_1(w_j)P_1(w_k)b=0$ by Lemma~\ref{lem:4.2},
since each $P_1(w_k)P_1(w_l)b$ must be positive. This proves that
$Q$ is faithful, and hence (a) holds.\\

(b)~ Let $B$ denote the $JBW^*$-algebra $P(Z)_2(v)$. Then by definition,\\
$B=\{\tpc{v}{\tp{v}{x}{v}_{P(Z)}}{v}_{P(Z)}:x\in P(Z)\}$. But
$$\tpc{v}{\tp{v}{x}{v}_{P(Z)}}{v}_{P(Z)}=P\tpc{v}{P\tp{v}{x}{v}}{v}=
P\tpc{v}{\tp{v}{x}{v}}{v}=PQ(v)^2x$$ so that $B=PQ(v)^2P(Z)$ and
\[
P_2(w)B=P_2(w)PQ(v)^2P(Z)=P_2(w)PP_2(w)Q(v)^2P(Z)=P_2(w)P(Z).
\]
Now let $F_v$ be the normal state space of
$B$, that is
\[
F_v=\{\ell\in B_*:\|\ell\|=1=\ell(v)\}.
\]
Recall from the first part of the proof that $v=Pw+PP_0(w)v$.
Also, $P(w)=P\tpc{v}{Pw}{v}=P(w)^\sharp$ implies that for $\ell\in
F_v$, $\ell(P(w))$ is real and
$1=\ell(v)=\ell(P(w))+\ell(PP_0(w)v)$. Therefore $\ell(P(w))\ge 0$
so that in fact $0\le P(w)\le v$, that is, $v-P(w)\in B^+$. Now
for each $i$, we have $f_i(v-Pw)=f_i(PP_0(w)v)=f_i(P_0(w)v)
=f_i(P_2(w_i)P_0(w)v)=0$. It follows, using Lemma~\ref{lem:4.2} as
above, that $v-Pw=0$.

Now for arbitrary $\ell\in F_v$, as $Pw=v$, we have
$\ell(w)=\ell(P(w))=\ell(v)=1$ and by \cite[Lemma 3.1]{FriRus83},
\begin{equation}\label{eq:neutral}
\ell=P_2(w)_*\ell.
\end{equation}
 By linearity and the
Jordan decomposition for self-adjoint functionals,
(\ref{eq:neutral}) extends to all $\ell\in B_*$. Hence for $b\in
B$, we have $\|b\|=\sup\{|\ell(b)|:\|\ell\|=1,\ell\in
B_*\}=\sup\{|\ell(P_2(w)b)|:\|\ell\|=1,\ell\in
B_*\}\le\|P_2(w)b\|$.  This proves (b). \qed

\begin{proposition}\label{prop:1prime}
Let $P$ be a normal contractive projection on a von Neumann algebra
$M$ of type I.  Then
$P(M)$ is a $JW^*$-triple of type I.
\end{proposition}
\pf\ Let $v$ be any nonzero tripotent of $P(M)$ and choose $w \in
M$ as in Lemma~\ref{lem:4.1}. By Proposition~\ref{prop:1},
$M_2(w)$ is a $JW^*$-algebra of Type I. By Lemma~\ref{lem:4.3} (a)
and \cite[Corollary 1.5]{EffSto79}, $P_2(w)P(M)=Q(M_2(w))$ is a
$JW^*$-subalgebra of $M_2(w)$, where $Q=P_2(w)P$ . By
Proposition~\ref{prop:3.1}, $Q(M_2(w))$ is a $JBW^*$-algebra of
type I, and by Lemma~\ref{lem:4.3}(b), $P(M)_2(v)$ (the Peirce
2-space of the tripotent $v$ of the $JW^*$-triple $P(M)$) is also
of type I since a unital surjective linear isometry is a Jordan
$^*$-isomorphism. One can now choose $v$ to be a complete
tripotent of $P(Z)$ to obtain from \cite[4.14]{H3} that $P(M)$ is
a $JW^*$-triple of type I.\qed

\section{Contractive projections on $JW^*$-triples}\label{sect:4}

\begin{proposition}\label{prop:3.3}
Let $Z$ be a $JBW^*$-triple of type I and let $v$ be a tripotent in
$Z$.  Then $P_2(v)Z$ is a $JBW^*$-algebra of type I.
\end{proposition}
\pf\ By Horn's structure theorem, we may assume that
$Z=L^\infty(\Omega,C)$ where $C$ is a Cartan factor. If $C$ is of
types 1,2, or 3, then there is a normal contractive projection $Q$
on $L^\infty(\Omega,\tilde{C})$, where $\tilde{C}$ is the von
Neumann envelope of $C$, with range $Z$. Since $P_2(v)Q$ is a
normal contractive projection  from the type I von Neumann algebra
$L^{\infty}(\Omega, \tilde{C})$ onto $P_2(v)Z$, the latter is of
type I by Proposition~\ref{prop:1prime}. If $C$ is of type 4, then
$P_2(v)Z=L^\infty(\Omega_2,C)\oplus L^\infty(\Omega_1)$, where
$\Omega_k=\{\omega\in\Omega: \mbox{ rank of }v(\omega)\mbox{ is }
k\}$, $k=0,1,2$. Indeed, if $f\in Z$ and $g=P_2(v)f$, then $g=0$
on $\Omega_0$, $g(\omega)=\langle
f(\omega),\widehat{v(\omega)}\rangle v(\omega)$ for
$\omega\in\Omega_1$ and $g=f$ on $\Omega_2$. Here we use the
notation $\hat{v}$ for the normal functional with support
tripotent $v$. It follows that the map $g=P_2(v)f \in P_2(v)Z
\mapsto (g_2,g_1)\in L^\infty(\Omega_2,C)\oplus
L^\infty(\Omega_1)$, where $g_1(\omega)= \langle
f(\omega),\widehat{v(\omega)}\rangle$ for $\omega\in\Omega_1$ and
$g_2=g|\Omega_2$, is a surjective linear isometry.

If $C$ is of types 5 or 6, then it is finite-dimensional and
$L^{\infty} (\Omega, C)$ is of type I$_{fin}$. By Lemma~\ref{sub},
the subtriple $P_2(v)(Z)$ is of the same type. \qed

\begin{theorem}\label{thm:2}  
Let $P$ be a normal contractive projection on a $JW^*$-triple $Z$ of
type I.  Then $P(Z)$ is of
type I.
\end{theorem}
\pf\
By \cite[4.14]{H3}, we need only show that $P(Z)_2(v)$ is of type I
for a complete tripotent $v \in P(Z)$.
Choose $w\in Z$ as in
Lemma~\ref{lem:4.1}.
By Proposition~\ref{prop:3.3}, $P_2(w)Z$ is a $JW^*$-algebra of type I.
One can now argue exactly as in the proof of
Proposition~\ref{prop:1prime},
 using
Lemma~\ref{lem:4.3}, to show that $P_2(w)P$ is a faithful, normal,
unital contractive projection of $P_2(w)Z$ onto $P_2(w)P(Z)$ (which is a
again subalgebra by \cite[Corollary 1.5]{EffSto79})
and that $P_2(w)$ is a unital isometry of $P(Z)_2(v)$ onto
$P_2(w)P(Z)$.  As in the proof of Proposition~\ref{prop:1prime}
and using Proposition~\ref{prop:3.1}, $P_2(w)P(Z)$ is of type I
and so is $P(Z)_2 (v)$. \qed

\begin{proposition}\label{prop:5.2}
Let $Z$ be a semifinite $JW^*$-triple and let $v$ be a partial
isometry in $Z$. Then $Z_2(v)$ is a semifinite $JW^*$-algebra.
\end{proposition}
\pf\
We prove this first in the case that $Z$ is a von Neumann algebra $M$.
If $M_2(v)$ had a type III part, we could follow $P_2(v)$ by the
projection
of $M_2(v)$ onto that type III part and obtain a Peirce 2-space of $M$
of type III. So we may assume that $M_2(v)$ is of type III. Let $p$ be
a finite nonzero projection in $M$ dominated by $v^*v$. Then $vp$ is a
nonzero projection
in $M_2(v)$ dominated by $v$ (cf. Proposition~\ref{prop:1}). We shall
show that $M_2(vp)$
is finite by showing that its involution is strongly continuous
on bounded spheres.

Let $x_\alpha$ be  a bounded net in $M_2(vp)$. Then
$$x_\alpha\stackrel{s}{\rightarrow} 0\mbox{ in }M_2(vp)
\Rightarrow vpx_\alpha^*vp(vp)^*x_\alpha\stackrel{w}{\rightarrow}0
\Rightarrow vpx_\alpha^* x_\alpha\stackrel{w}{\rightarrow}0
\Rightarrow $$ $$ px_\alpha^* x_\alpha\stackrel{w}{\rightarrow}0
\Rightarrow px_\alpha^* x_\alpha p\stackrel{w}{\rightarrow}0
\Rightarrow x_\alpha p\stackrel{s}{\rightarrow}0
\Rightarrow(\mbox{ by \cite[p. 97-98]{SA}}) $$  $$
x_\alpha^*=px_\alpha^* \stackrel{s}{\rightarrow}0 \Rightarrow
x_\alpha x_\alpha^* \stackrel{w}{\rightarrow}0 \Rightarrow
x_\alpha x_\alpha^*vp \stackrel{w}{\rightarrow}0 \Rightarrow
x_\alpha (vp)^*vpx_\alpha^*vp \stackrel{w}{\rightarrow}0 $$  $$
\Rightarrow x_\alpha \circ x_\alpha^\sharp
\stackrel{w}{\rightarrow}0 \Rightarrow x_\alpha^\sharp
\stackrel{s}{\rightarrow}0\mbox{ in }M_2(vp). $$ Thus $vp$ is a
finite projection which is a contradiction.

To prove the general case,
write $Z=Z_I\oplus Z_{II}$ where $Z_I$ is of type I and $Z_{II}$
is of type II. Since $P_2(v)Z=P_2(v_1)Z_I\oplus P_2(v_2)Z_{II}$ for
suitable partial isometries $v_1\in Z_I$ and $v_2\in Z_{II}$, and we
already know that $P_2(v_1)Z_I$ is of type I, we may assume by
\cite{H2} that $Z$ is triple isomorphic to $pM\oplus H(N,\alpha)$, where
$M$ and $N$ are
von Neumann algebras of type II. Accordingly $v_2=v_2'+v_2''$
so that $P_2(v_2)Z_{II}=P_2(v_2')(pM)\oplus P_2(v_2'')(H(N,\alpha))=
M_2(v_2'')\oplus H(N_2(v_2''),\alpha)$ and
$Q(N_2(v_2''))=H(N_2(v_2''),\alpha)$
where $Q$ is the projection $Q(x)=(x+\alpha(x))/2$ for $x\in N$.
By the first part of the proof, both $M_2(v_2')$ and $N_2(v_2'')$
are semifinite. Then by Proposition~\ref{prop:3}, $P_2(v_2'')H(N,\alpha)$
is semifinite and the result follows.\qed

\begin{theorem}\label{thm:3}
Let $P$ be a normal contractive projection on a semifinite $JW^*$-triple
$Z$. Then $P(Z)$ is a semifinite $JW^*$-triple.\end{theorem}
\pf\
By passing to the type III part of $P(Z)$,
assuming it is nonzero for contradiction,
 and using \cite{H2},
we may assume that $P(Z)=pM\oplus H(N,\alpha)$ where $M$ and $N$ are
von Neumann algebras of type III.
As in the proof of Proposition~\ref{prop:1prime},
using Lemma~\ref{lem:4.3}, Proposition~\ref{prop:5.2},
and Proposition~\ref{prop:3}, one shows that
 $P(Z)_2(v)$ is semifinite for any tripotent $v$
of $P(Z)$.   Choosing $v=0\oplus 1_N$ leads to
$P(Z)_2(v)=H(N,\alpha)$, a contradiction unless
$H(N,\alpha)=0$.  Choosing $v=p\oplus 0$ leads
to $P(Z)_2(v)=pMp$, again a contradiction unless $pMp=0$,
which implies that $p=0$, another contradiction.
\qed

\section{Contractive projections on $JBW^*$-triples}\label{sect:5}

In this section we extend Theorems~\ref{thm:2} and ~\ref{thm:3} to
arbitrary $JBW^*$-triples and make some remarks
on the atomic case.

A close examination of the proof of Theorem~\ref{thm:2}
reveals that it carries over to the case of $JBW^*$-triples if we can show
that the range of a faithful normal positive unital
projection on a Type I JBW*-algebra is Type I.  For this,
one only needs to decompose a $JBW$-algebra
$A$ into a direct sum of a $JC$-algebra $A_{sp}$ and an ``exceptional
algebra''
$A_{ex}$ of the form
$L^\infty(\Omega,M_3^8)$ where $M_3^8$ denotes the $JBW$-algebra whose
complexification
is the Cartan factor of type 6 (\cite{Shultz79}). Since both of these
summands are
unital subalgebras, it follows easily that
 any positive unital projection on $A$ restricts to a positive unital
projection on each summand. The image of the restriction of a
faithful projection to $A_{sp}$ is a subalgebra of $A_{sp}$ by
\cite[Corollary 1.5]{EffSto79}, and we can apply
Proposition~\ref{prop:3.1}.  The image of the restriction of the
projection to $A_{ex}$ is of type I by the following lemma which
applies {\it verbatim} to type I$_{fin}$ JBW-algebras of which
$A_{ex}$ is one.

\begin{lemma}\label{dpp}
Let $Z$ be a type I$_{fin}$ JBW*-triple and let $P: Z \longrightarrow
Z$ be a normal contractive projection. Then $P(Z)$ is a type I$_{fin}$
JBW*-triple.
\end{lemma}
\pf\
We note that $P(Z)$ is norm-closed. By weak* continuity of $P$ and the
Krein-Smulyan Theorem, $P(Z)$ is also  weak* closed.
Also, $P$ induces a contractive projection $P_* : f \in Z_*
\mapsto f \circ P \in Z_*$ on the predual $Z_*$. By \cite{CM},
$Z_*$ has the Dunford-Pettis property. The predual of $P(Z)$ identifies
with $Z_* /P_*^{-1}(0)$ which is linearly isometric to the complemented
subspace $P_*(Z_*)$ of $Z_*$, and therefore has the Dunford-Pettis
property. Hence by \cite{CM} again, $P(Z)$ is of type I$_{fin}$.
\qed

Now, proceeding exactly as in the proof of Theorem~\ref{thm:2}  we have
the result for JBW*-triples.
\begin{theorem}
Let $P$ be a normal contractive projection on a $JBW^*$-triple $Z$ of
type I.  Then $P(Z)$ is of
type I.
\end{theorem}

As noted before, a type II JBW*-triple is a JW*-triple. It follows from
Proposition~\ref{prop:3.3} and Theorem~\ref{thm:3} that if $Z$ is
a semifinite JBW*-triple,
then $P_2(e)Z$ is also a semifinite JBW*-algebra. Using this fact,
now there is no difficulty of extending the proof of Theorem~\ref{thm:3}
to the case of JBW*-triples.

\begin{theorem}
Let $P$ be a normal contractive projection on a semifinite $JBW^*$-triple
$Z$. Then $P(Z)$ is a semifinite $JBW^*$-triple.\end{theorem}

Tomiyama \cite{TO} has proved that the a von Neumann algebra which
is the range of a normal contractive projection on an atomic von
Neumann algebra is itself atomic.  It is also known (see
\cite[Exercise 8, p.334]{Tak79}) that a von Neumann algebra
$M\subset B(H)$ is atomic if  and only if there is a faithful
family of normal conditional expectations of $B(H)$ onto $M$. We
end with a very simple proof of the following result which extends
Tomiyama's theorem to $JBW^*$-triples. The proof follows from a
result in \cite{CI} which states that a JBW*-triple is atomic if,
and only if, its predual has the Radon-Nikodym property. The
following result is clearly false without the normality assumption
on $P$.

\begin{proposition}\label{atomic} Let $Z$ be an atomic
JBW*-triple and let
$P: Z\longrightarrow Z$ be a normal contractive
projection. Then the range $P(Z)$ is (linearly isometric to) an
atomic JBW*-triple.
\end{proposition}
\pf\ As in the proof of Lemma~\ref{dpp},
 the predual of $P(Z)$
is linearly isometric to a complemented subspace of the predual $Z_*$
which has the Radon-Nikodym property. So $P(Z)$ is atomic by \cite{CI}.
\qed

\it
\begin{tabbing}

 Goldsmiths College1111111111111111     \=     University of
California\kill
 C-H. Chu \> M. Neal and B. Russo\\
 Goldsmiths College \> University of California\\
University of London \>Irvine, California \\ London SE14 6NW \>
92697-3875\\U.\ K. \> U.\ S.\ A.
\end{tabbing}

\end{document}